\documentstyle{article}
\input amssym.def

\newcommand{\N} {{\Bbb N}}
\newcommand{\Q} {{\Bbb Q}}
\newcommand{\R} {{\Bbb R}}

\newcommand{\Z} {{\Bbb Z}}
\newcommand{\bi}{\bibitem}

\newcommand{\mathematica}{Mathematica}

\newcommand{\reduce}{REDUCE}
\newcommand{\1}{{\bf{1}}}
\newcommand{\2}{{\bf{2}}}
\newcommand{\4}{{\bf{4}}}
\newcommand{\5}{{\bf{5}}}
\newcommand{\6}{{\bf{6}}}
\newcommand{\7}{{\bf{7}}}
\newcommand{\8}{{\bf{8}}}

\newcommand{\0}{{\bf{0}}}

\begin{document}
\noindent
{\Large \bf The Algebra of Holonomic Equations}
\\[3mm]
{\bf
Wolfram Koepf
}
\\[1mm]
Konrad-Zuse-Zentrum Berlin, Heilbronner Str.\ 10, D--10711 Berlin
\\[5mm]
%
{\bf Abstract.}
In this article algorithmic methods are presented that have essentially
been introduced into computer algebra systems like Maple or Mathematica
within the last decade. The main
ideas are due to Stanley and Zeilberger.
Some of them had already been discovered in the last century by Beke,
but because of their complexity 
the underlying algorithms have fallen into oblivion.
We give a survey of these techniques, show how they can be used to
identify transcendental functions, and
present implementations of these algorithms in
computer algebra systems.

\section{Algebraic Representation of Transcendental Functions}

How can transcendental function be represented by algebraic means?
To give this question another flavor:
What is the main difference between the exponential function $f(x)=e^x$
and the function $g(x)=e^x+|x|/10^{1000}$, that makes $f$ an
elementary function, but not $g$, although $f$ and $g$ are numerically
quite close on a part of the real axis?

Or let's consider
an example of discrete mathematics: Why is the factorial function
$a_n=n!$ considered to be the most important discrete function, and not 
$b_n=n!+n/10^{1000}$ or any other discrete function?

Although these examples refer to the most important continuous and discrete
{\sl transcendental} functions, oddly enough the answers to the above
questions are {\sl purely algebraic}: The exponential function $f$
is characterized by any of the following algebraic properties:
\begin{enumerate}
\item
$f$ is continuous, $f(1)=e$, and for all $x,y$ we have $f(x+y)=f(x)\cdot f(y)$;
\item
$f$ is differentiable, $f'(x)=f(x)$ and $f(0)=1$;
\item
$f\in C^{\infty}$, $f(x)=\sum\limits_{n=0}^{\infty}a_n\,x^n$ with
$a_0=1$, and for all $n\geq 0$ we have $(n+1)\,a_{n+1}=a_n$;
\end{enumerate}
and the factorial function $a_n$
is represented by any of the following algebraic properties:
\begin{enumerate}
\item[4.]
$a_0=1$, and for all $n\geq 0$ we have $a_{n+1}=(n+1)\,a_n$;
\item[5.]
the generating function $f(x)\!=\!\sum\limits_{n=0}^{\infty}a_n x^n$
satisfies the differential equation \newline
$x^2 f'(x)+(x-1)f(x)+1\!=\!0$
with the initial condition $f(0)=1$. 
\end{enumerate}
(Note here that one could argue
that property (1.) is not algebraic since the symbol $e$ is needed
in the representation.) I do not know any method to represent
transcendental functions using functional equations, such as property (1.),
but 
I will show, why and how the other properties can be suitable for this purpose, 
being mainly concerned with properties (2.) and (4.).
In \S~\ref{sec:Holonomic Systems of Several Variables} we consider, how these 
representations can be viewed as purely polynomial cases.

Observe that the ``generating function'' of the factorial function is
convergent only at the origin, and therefore must be considered as a
formal series. In particular, a ``closed representation'' (whatever
that should mean) of the generating function cannot be given. 
But this is not the main issue here. Rather than working with the
generating function itself, it is much better to work with
its differential equation which is purely algebraic (in fact, it is
purely polynomial). The same argument applies to the exponential
and factorial functions themselves. Rather than working with these
transcendental objects, one should represent them by their corresponding
differential and recurrence equations.

The given properties are {\sl structural statements} about the
corresponding functions. Any small modification (even changing the
value at a single point) destroys this structure.
For example, the function $g(x)=e^x+|x|/10^{1000}$ {\sl cannot} be
characterized by a rule analogous to one of the properties (2.)--(3.).
On the other hand, the function $h(x)=e^x+x/10^{1000}$ can be represented
by the differential equation $(x-1)\,h''(x)-x\,h'(x)+h(x)=0$ with the initial
values $h(0)=1$ and $h'(0)=1+10^{-1000}$.

Therefore, the special (and common) fact about the exponential and factorial 
functions is that they both satisfy a differential or recurrence equation,
respectively, that is homogeneous, linear, of order one, and has polynomial 
coefficients.

We can generalize this observation \cite{Zei1}:
A continuous function of one variable $f(x)$ is {\sl holonomic}, 
if it satisfies a homogeneous linear
differential equation with polynomial coefficients; 
we call such a differential equation also holonomic.

By linear algebra arguments,
Stanley \cite{Sta} showed that sums and products of holonomic functions
and the composition with algebraic functions also form holonomic functions.
This can be seen as follows: Assume $f$ and $g$ satisfy holonomic
differential equations of order $n$ and $m$, respectively.
We consider the linear space $L_f$ of functions with rational coefficients
generated by $f, f',f'',\ldots,f^{(k)},\ldots$. Since $f,
f',\ldots,f^{(n)}$ are linearly dependent by the given
holonomic differential equation
and since by differentiation the same conclusion follows for 
$f',f'',\ldots,f^{(n+1)}$, and so on inductively,
the dimension of $L_f$ is $\leq n$. Similarly $L_g$ has dimension $\leq m$.
We now build the sum $L_f +L_g$ which is of dimension $\leq n+m$.
As $f\!+\!g,(f\!+g)',\ldots,(f\!+\!g)^{(k)},\ldots$ are
elements of $L_f +L_g$, arbitrary $n+m+1$ many of them are linearly dependent.
In particular, $f+g$ satisfies a holonomic differential
equation of order $\leq n+m$.

Similarly the product and composition cases can be handled.
Note that the above proof provides a construction of the resulting holonomic
equation by linear algebra techniques.
It is remarkable that 100 years ago, Beke \cite{Beke1}--\cite{Beke2} 
already described these
algorithms to generate holonomic differential equations for the sum and 
product of $f$ and $g$ from the holonomic differential equations of $f$ and 
$g$. Hence, he had discovered algorithmic versions of Stanley's results!

Analogously, a discrete function (sequence) of one variable
is called holonomic, if it satisfies a
homogeneous linear recurrence equation with polynomial coefficients.
Such a recurrence equation is also called holonomic.
Sums and products of discrete holonomic functions are again holonomic,
and there are similar algorithms to calculate representing holonomic
recurrence equations (s.\ \cite{SZ}, \cite{KS}).

A function 
\[
f(x)=\sum\limits_{n=0}^{\infty}a_n\,x^n
\]
represented by a power series is holonomic if and only if 
the corresponding power series coefficient $a_n$ is a holonomic sequence.
The holonomic equations for $f(x)$ and $a_n$ can be converted 
equating coefficients.

Note that these algorithms were implemented by Salvy and Zimmermann
in the {\tt gfun} package of Maple's share library \cite{SZ}. I wrote
a Mathematica implementation, {\tt SpecialFunctions},
to be obtained by World Wide Web from the address
{\tt ftp://ftp.zib-berlin.de/pub/UserHome/Koepf/SpecialFunctions}.
\linebreak
Examples of this implementation will be given later.

\section{Identification of Transcendental Functions}
\label{sec:Identification of Transcendental Functions}

Note that the notion of holonomy provides a {\sl normal form}
for a suitably large number of transcendental functions,
which can then be utilized for {\sl identification} purposes.
The holonomic equation {\sl of lowest order} corresponding to a
holonomic function constitutes such a normal form. Once we have
calculated the normal form of a holonomic function, the latter is identified:
Two holonomic functions are identical if and only if they have the
same normal form, and satisfy the same initial conditions.

But also without having access to the {\sl lowest order} holonomic equations, 
one can check whether two holonomic functions agree, since (by linear algebra,
e.g.,) it is easy to see whether two holonomic equations are compatible
with each other.

Therefore, we may ignore that $e^x, \sin x, \cos x, \arctan x, \arcsin x$ 
and others form transcendental functions, and take only their
holonomic differential equations 
$f'=f$, $f''=-f$, $f''=-f$, $(1+x^2)f''+2xf'=0$, $(x^2-1)f''+xf'=0$ etc.
into account. From these differential equations, corresponding 
differential equations for sums and products can be generated by 
the above mentioned technique, using only polynomial arithmetic and 
linear algebra.
For example, the function $f(x)=\arcsin^2 x$ yields
$(x^2-1)f'''+3xf''+f'=0$. Note, however,
that in the given case one can get even more: 
The resulting holonomic differential equation 
is directly equivalent to the holonomic recurrence equation
$n(1 + n) (2 + n)a_{n+2}=n^3a_n$ for the coefficients
$a_n$ of the Taylor series of  $\arcsin^2 x=\sum\limits_{n=0}^\infty a_n x^n$,
and since this holonomic recurrence equation fortunately 
contains only the two terms
$a_{n+2}$ and $a_n$, it can be solved explicitly, and leads to the
representation
\[
\arcsin^2 x=\sum_{n=0}^\infty
\frac{4^n\,n!^2}{(1 + n)\,(1 + 2n)!}x^{2n+2}
\]
(compare \cite{GKP}, \cite{Wil1}, \cite{Koe92}--\cite{Koe93}).

Note that not only a function like the Airy function ${\rm Ai}\:(x)$
(s.\ e.\ \cite{AS}, (10.4)) 
falls under the category of holonomic functions, since it satisfies
the simple holonomic differential equation $f''-xf=0$, moreover 
the classical families of orthogonal polynomials%
\footnote{As families of orthogonal polynomials they are {\sl not}
polynomials!}
and many other special functions form holonomic
functions \cite{AS}. These depend on several variables, and we will discuss 
this situation in \S~\ref{sec:Holonomic Systems of Several Variables}.

On the other hand, there are functions that are not holonomic, like
the tangent function $\tan x$ (s.\ \cite{Sta}, \cite{KS}).
The identification problem for
expressions involving nonholonomic functions can only be treated
after preprocessing the input. If, for example, we want to verify the addition 
formula for the tangent function 
\[
\tan\:(x+y)={{\tan x + \tan y}\over {1 - \tan x\,\tan y}}
\]
by the given method, then we have to replace all occurrences
of the tangent function by sines and cosines (which are holonomic)
using the rewrite rule $\tan x=\sin x/\cos x$.
We can then generate a polynomial equation by multiplying both sides
by the common denominator.
This procedure results in the equivalent representation
\begin{equation}
\left( \cos x\cos y - \sin x\sin y \right) \sin\: (x + y)
=
\left( \cos y\sin x + \cos x\sin y \right) \cos \:(x + y)
\label{eq:tanadditionformula}
\end{equation}
which is easily proved since the algorithms generate the common holonomic 
differential equation $f''(x)+4f'(x)= 0$ with respect to $x$ (or the
common holonomic differential equation $f''(y)+4f'(y)= 0$ with respect to $y$)
for both sides of (\ref{eq:tanadditionformula}) where the common initial values
are $f(0)=\cos y\,\sin y$, and $f'(0)={{\cos y}^2} - {{\sin y}^2}$.
Assume that for the initial value functions we had obtained
different representations
(e.g.\ $\cos y\,\sin y$ and $\sin\:(2y)/2$). These could be verified
by the same technique.

In the \mathematica\ package {\tt SpecialFunctions} 
(s.\ also \cite{KoepfGAMM}), the procedure
{\tt HolonomicDE[f,x]} calculates the holonomic differential equation
of $f$ with respect to the variable $x$ using the known holonomic differential 
equations of the primitive functions, and the sum and product algorithms
by recursive decent through the expression tree. 
Here we call a function {\sl primitive} if it is rational, or
whenever we use a separate symbol for it
and a holonomic differential equation is known. Therefore the above 
mentioned functions (besides the tangent function) are primitive.

The examples given are governed by the following \mathematica\ session:

{\small
\begin{verbatim}
In[1]:= <<SpecialFunctions`

In[2]:= HolonomicDE[ArcSin[x]^2,x]

                                               (3)
Out[2]= F'[x] + 3 x F''[x] + (-1 + x) (1 + x) F   [x] == 0

In[3]:= DEtoRE[%,F,x,a,n]

         3
Out[3]= n  a[n] - n (1 + n) (2 + n) a[2 + n] == 0

In[4]:= Series[ArcSin[x]^2,{x,0}]

              k  2 + 2 k   2
             4  x        k!
Out[4]= Sum[------------------, {k, 0, Infinity}]
            (1 + k) (1 + 2 k)!

In[5]:= HolonomicDE[AiryAi[x],x]

Out[5]= -(x F[x]) + F''[x] == 0

In[6]:= HolonomicDE[AiryAi[x]^2,x]

                              (3)
Out[6]= 2 F[x] + 4 x F'[x] - F   [x] == 0

In[7]:= HolonomicDE[Sin[x+y]*(Sin[x]Sin[y]-Cos[x]Cos[y]),x]     

                   (3)
Out[7]= 4 F'[x] + F   [x] == 0

In[8]:= HolonomicDE[Cos[x+y]*(Sin[x]Cos[y]+Cos[x]Sin[y]),x]

                   (3)
Out[8]= 4 F'[x] + F   [x] == 0

In[9]:= HolonomicDE[Cos[y]*Sin[y],y]

                   (3)
Out[9]= 4 F'[y] + F   [y] == 0

In[10]:= HolonomicDE[Sin[2y]/2,y]

Out[10]= 4 F[y] + F''[y] == 0
\end{verbatim}
}\noindent
One difficulty that may arise with the method described is 
that in some instances
the sum and product algorithms will not generate the holonomic
differential equation of lowest order, as in the above example for 
$\cos y\sin y$. In this case, the normal form property is lost.
In fact, the sum algorithm calculates a holonomic equation that is valid
for {\sl any linear combination} $af+bg$ rather than the particular given sum
$f+g$.
As a simple example, we consider the sum $\sqrt {1+x}+{\frac {1}{\sqrt {1+x}}}$
satisfying the first order differential equation
\[
2\,\left (2+x\right )\left (1+x\right )F'(x)-xF(x)=0
\;.
\]
This differential equation can be found using a method given in
\cite{Koe92}--\cite{Koe93}, whereas the sum algorithm generates the
second order differential equation
\[
4\,{{\left( 1 + x \right) }^2}\,F''(x)+
4\,\left( 1 + x \right) \,F'(x)-F(x)= 0
\;.
\]
The reason for the existence of a differential equation
of lower order is due to the fact that the
ratio of the two summands $\sqrt {1+x}$ and ${\frac {1}{\sqrt {1+x}}}$
forms a rational function.

Similarly, the sum of two consecutive Legendre polynomials
$P_n(x)+P_{n+1}(x)$ satisfies the second order differential equation
\[
\left (x-1\right )\left (x+1\right )F''(x)+
\left (x+1\right )F'(x)-\left (n+1\right )^{2}F(x)=0
\;,
\]
whereas the sum algorithm generates the differential equation
\begin{eqnarray*}
0&=&
{{\left( x-1 \right) }^2}\,{{\left( 1 + x \right) }^2}\,F''''(x) +
8\,\left( x-1 \right) \,x\, \left( 1 + x \right) \,F'''(x)
\\&&
+
2\,\left( -2 + 2\,n + {n^2} + 6\,{x^2} - 2\,n\,{x^2} - {n^2}\,{x^2}
        \right) \,F''(x) 
\\&&
- 4\,n\,\left( 2 + n \right) \,x\,F'(x) + 
n\,{{\left( 1 + n \right) }^2}\,\left( 2 + n \right) \,F(x)
\end{eqnarray*}
of fourth order, which is also valid 
for the difference $P_n(x)-P_{n+1}(x)$ and for any
other linear combination.

For the {\sl verification of identities}, this is not an important issue,
since the {\sl compatibility} of two holonomic equations can be easily checked.
This situation is similar to
proving a rational identity by pure polynomial arithmetic
without gcd computations (after having multiplied through by all denominators),
and is actually equivalent to a noncommutative polynomial division, see 
\S~\ref{sec:Holonomic Systems of Several Variables}.

In the case that the normal form is needed for a particular problem,
a {\sl factorization algorithm} 
can be used, s.\ 
\S~\ref{sec:Noncommutative Factorization and Holonomic Normal Form}. 

For the discrete functions, the situation is quite similar. We
call a function primitive whenever we use a separate symbol for it
and a holonomic recurrence equation is known. To these primitive functions,
we add the rational functions and the functions
\begin{equation}
(mn+b)!\;,
\quad
\frac{1}{(mn+b)!}\;
\quad(m\in\Q)
\;,
\quad
\mbox{and}
\quad
a^n
\label{eq:zulaessige diskrete Funktionen}
\end{equation}
whose holonomic recurrence equations are known,
as primitive functions with respect to the variable $n$.
We consider the factorial function to be equivalent to the $\Gamma$ function
$\Gamma\:(a+1)=a!$, and declare binomial coefficients etc.\ also via
factorials.
 From the holonomic representations of the primitive functions the
holonomic equations for all sums and products can be established.
E.\ g.\ the two equations
\begin{equation}
(n-k+1)^2 F(n+1,k)-(1+n)^2 F(n,k)=0
\label{eq:combnk^2.1}
\end{equation}
and
\begin{equation}
(k+1)^2 F(n,k+1)-(n-k)^2 F(n,k)=0
\label{eq:combnk^2}
\end{equation}
for $F(n,k)={{n}\choose{k}}^2$. Whereas these are simple consequences
of the representation of $F(n,k)$ by factorials, the given procedure
can be applied, for example, to the more complicated function 
$F(n,k)=\frac{n!+k!^2}{k}$ to generate the two holonomic equations
\[
n F(n+2,k) -(1+3n+n^2) F(n+1,k) +(1+n)^2 F(n,k)=0
\]
and
\[
k (2 + k)^2 F(n,k+2)
-(1 + k) (1 + 3 k + k^2)
(3 + 3 k + k^2) F(n,k+1)
+ k (1 + k)^3 F(n,k)
=0
\;.
\]
Note that the given approach also covers
all kinds of orthogonal polynomials and special functions with respect
to their discrete variables, see 
\S~\ref{sec:Holonomic Systems of Several Variables}.

In our \mathematica\ implementation {\tt SpecialFunctions}, the procedure
\linebreak
{\tt HolonomicRE[a,n]} calculates the holonomic recurrence equation
of $a_n$ with respect to the variable $n$ taking the known holonomic recurrence
equations of the primitive functions into account, 
and using the sum and product algorithms
by recursive decent through the expression tree. The above examples are
governed by the following \mathematica\ session:

{\small
\begin{verbatim}
In[11]:= HolonomicRE[Binomial[n,k]^2,n]

                2                   2
Out[11]= (1 + n)  a[n] - (1 - k + n)  a[1 + n] == 0

In[12]:= HolonomicRE[Binomial[n,k]^2,k]

                 2               2
Out[12]= (-k + n)  a[k] - (1 + k)  a[1 + k] == 0

In[13]:= HolonomicRE[(n!+k!^2)/k,n]

                2                     2
Out[13]= (1 + n)  a[n] + (-1 - 3 n - n ) a[1 + n] + n a[2 + n] == 0

In[14]:= HolonomicRE[(n!+k!^2)/k,k]

                  3
Out[14]= k (1 + k)  (3 + k) a[k] - 
 
                        2              2 
>   (1 + k) (1 + 3 k + k ) (3 + 3 k + k ) a[1 + k] + 

             2
>   k (2 + k)  a[2 + k] == 0
\end{verbatim}
}\noindent

\section{Hypergeometric Sums}

Rather than having functions given as finite sums and products of
primitive expressions, an important
class of functions, particularly in analysis and 
combinatorics, is given by infinite sums of products 
of terms of the form (\ref{eq:zulaessige diskrete Funktionen})
\begin{equation}
s(n)=\sum_{k\in\Z} F(n,k)
\;.
\label{eq:infinitesum}
\end{equation}
Then $F(n,k)$ is an $(m,l)$-fold {\sl hypergeometric term}. That is, both
$F(n+m,k)/F(n,k)$ and $F(n,k+l)/F(n,k)$ are rational functions with
respect to $n$ and $k$ for a certain pair $(m,l)\in\N^2$.
For example, by (\ref{eq:combnk^2.1})--(\ref{eq:combnk^2})
this is valid for $F(n,k)={{n}\choose{k}}^2$ with $m=l=1$.
We assume moreover that the sums (\ref{eq:infinitesum}) have finite support,
i.e., they are finite sums for each particular $n\in\N$.

A modification \cite{Koe94z} of
the (fast) {\sl Zeilberger algorithm}
(\cite{Zei2}, see also \cite{Koornwinder}, and \cite{PS})
returns a holonomic recurrence equation valid for $s(n)$.
Zeilberger's algorithm is based on a decision procedure for
indefinite summation due to Gosper \cite{Gos}. In our example case, Zeilberger's
algorithm finds the holonomic recurrence equation
$(1+n)\,s(n+1)=2(1+2n)\,s(n)$ for 
$s(n)=\sum\limits_{k\in\Z} {{n}\choose{k}}^2=\sum\limits_{k=0}^n {{n}\choose{k}}^2$ 
which fortunately has only two terms. Therefore, we are led to the
representation
\[
s(n)=\sum\limits_{k=0}^n {{n}\choose{k}}^2=\frac{(2n)!}{n!^2}
\;.
\]
Even though, in general, the resulting recurrence 
equation 
has more than two terms, this holonomic equation
contains very important structural information
about $s(n)$. This may be used to show that a certain family of polynomials
is orthogonal or not \cite{Zei3}, and can be an interesting property for 
numerical purposes (compare \cite{Deufl1}--\cite{Deufl2}).

In particular, as described in the last section,
the generated structural information can be used for the
identification of a transcendental function that is
given as sum (\ref{eq:infinitesum}). Note that sums of type
(\ref{eq:infinitesum}) in general form transcendental functions with
respect to the discrete variable $n$.

For example, to check the identity (compare \cite{Strehl2})
\begin{equation}
\sum_{k=0}^n
{{n}\choose{k}}^3
=
\sum_{k=0}^n
{{n}\choose{k}}^2 {{2k}\choose{n}}
\label{eq:disguise}
\end{equation}
which is nontrivial since for $n=1$ it reads $1+1=0+2$,
we need only to show
that both sums $s(n)$ satisfy the common recurrence equation
\begin{equation}
(n+2)^2 s(n+2)-(16 + 21 n + 7 n^2)s(n+1)-(n+1)^2 s(n)=0
\label{eq:disguise2}
\end{equation}
which is the result given by Zeilberger's algorithm. We also have
the same initial values $s(0)=1$ and $s(1)=2$, so we are done.

In Mathematica these computations are done by

{\small
\begin{verbatim}
In[15]:= HolonomicRE[Sum[Binomial[n,k]^2,{k,0,n}],n]

Out[15]= -2 (1 + 2 n) a[n] + (1 + n) a[1 + n] == 0

In[16]:= HolonomicRE[Sum[Binomial[n,k]^3,{k,0,n}],n]

                   2                         2
Out[16]= -8 (1 + n)  a[n] + (-16 - 21 n - 7 n ) a[1 + n] + 
 
             2
>     (2 + n)  a[2 + n] == 0

In[17]:= HolonomicRE[Sum[Binomial[n,k]^2*Binomial[2k,n],{k,0,n}],n]

                   2                         2
Out[17]= -8 (1 + n)  a[n] + (-16 - 21 n - 7 n ) a[1 + n] + 
 
             2
>     (2 + n)  a[2 + n] == 0
\end{verbatim}
}\noindent
Note that the example shows that transcendental functions can come in
quite different disguises. Might the left or the right hand side
of (\ref{eq:disguise}) be a preferable representation? 
This question cannot be answered satisfyingly. 
A holonomic recurrence equation like
(\ref{eq:disguise2}), defining the same transcendental function $s(n)$,
is probably the simplest way to describe a function of a discrete variable,
since it postulates how the values of the function can be calculated
iteratively. Not only is this a quite efficient way to calculate the values
of $s(n)$, but moreover it is preferable to either of the two
representations given in (\ref{eq:disguise}), since it gives a
unique representation scheme. This is what a normal form is about.

As a further example, we consider the function
($\alpha,\beta,\gamma\in\N_0,\;z,M,d\in\R^+$)
\begin{eqnarray*}
V(\alpha,\beta,\gamma)
&=&
(-1)^{\alpha+\beta+\gamma}\cdot
\frac{\Gamma(\alpha\!+\!\beta\!+\!\gamma\!-\!d)
\Gamma(d/2\!-\!\gamma)\Gamma(\alpha\!+\!\gamma\!-\!d/2)
\Gamma(\beta\!+\!\gamma\!-\!d/2)}
{\Gamma(\alpha)\Gamma(\beta)\Gamma(d/2)\Gamma(\alpha+\beta+2\gamma-d)
M^{\alpha+\beta+\gamma-d}}
\\\\&&\cdot
\; _2 F_1\left.
\!\!
\left(
\!\!\!\!
\begin{array}{c}
\multicolumn{1}{c}{\begin{array}{cc} \alpha+\beta+\gamma-d
\;, & \alpha+\gamma-d/2 \end{array}}\\[1mm]
\multicolumn{1}{c}{ \alpha+\beta+2\gamma-d}
            \end{array}
\!\!\!\!
\right| z\right)
\end{eqnarray*}
($ _2 F_1$ here represents Gau{\ss}'s hypergeometric function,
see \cite{AS}, Chapter 15), which plays a role for the 
computation of Feynman-diagrams \cite{FT}%
\footnote{I am indebted to Jochem Fleischer who informed me about
a misprint in formula (31) of \cite{FT}.}, for which Zeilberger's
algorithm generates the holonomic recurrence equation
\begin{eqnarray*}
0&=&
  \left( \alpha  + \beta  - d + \gamma  \right) \,
     \left( 2\,\alpha  - d + 2\,\gamma  \right) \,V(\alpha,\beta,\gamma ) 
\\&&
+\; 
    \alpha \,M\,\left( 2\,\alpha  + 2\,\beta  - 2\,d + 4\,\gamma  - 2\,z - 
       4\,\alpha \,z - 2\,\beta \,z + 3\,d\,z - 4\,\gamma \,z \right) \,
V(\alpha+1,\beta,\gamma ) 
\\&&
+\; 2\,\alpha \,\left( 1 + \alpha  \right) \,{M^2}\,
     \left( z-1 \right) \,z\,V(\alpha+2,\beta,\gamma ) 
\end{eqnarray*}
and analogous recurrence equations with respect to the variables $\beta$ 
and $\gamma$ (see \cite{Koepf95}). 
These, in particular, can be used for numerical purposes.

Note that for the application of Zeilberger's algorithm our 
\mathematica\ program uses the Paule-Schorn implementation \cite{PS}.
For the current example, the output is given by

{\small
\begin{verbatim}
In[18]:= HolonomicRE[(-1)^(alpha+beta+gamma)*Gamma[alpha+beta+gamma-d]*
         Gamma[d/2-gamma]*Gamma[alpha+gamma-d/2]*Gamma[beta+gamma-d/2]/
         (Gamma[alpha]*Gamma[beta]*Gamma[d/2]*
         Gamma[alpha+beta+2*gamma-d]*M^(alpha+beta+gamma-d))*
         Hypergeometric2F1[alpha+beta+gamma-d,alpha+gamma-d/2,
         alpha+beta+2*gamma-d,z],alpha,V]

Out[18]= (alpha + beta - d + gamma) (2 alpha - d + 2 gamma) V[alpha] + 
 
>     alpha M (2 alpha + 2 beta - 2 d + 4 gamma - 2 z - 4 alpha z - 
 
>        2 beta z + 3 d z - 4 gamma z) V[1 + alpha] + 
 
                           2
>     2 alpha (1 + alpha) M  (-1 + z) z V[2 + alpha] == 0
\end{verbatim}
}\noindent

\section{Holonomic Systems of Several Variables}
\label{sec:Holonomic Systems of Several Variables}

In \cite{Zei1}, Zeilberger considered the more general situation of
functions $F$ of several discrete and continuous variables. If we have
$d$ variables, and $d$ (essentially independent) mixed homogeneous linear
(partial) difference-differential equations with polynomial coefficients
in all variables are given for $F$, then $F$ is called a {\sl holonomic system}
(compare \cite{Bernstein1}--\cite{Bj}). In most cases 
these holonomic equations
together with suitably many initial values declare $F$ uniquely.

In particular, we concentrate on the situation, when the given system
of holonomic equations is separated, i.e.\ each of them is either an
ordinary differential equation or a pure recurrence equation. These
representing holonomic equations can be generated by the method
described in
\S~\ref{sec:Identification of Transcendental Functions}
whenever $F$ is given in terms of sums and products of primitive functions.

For example,
the Legendre polynomials $F(n,x)=P_n(x)$ (\cite{AS}, Chapter 22) form a
holonomic system by their holonomic differential equation
\begin{equation}
(x^2-1)F''(n,x)+2xF'(n,x)-n(1+n)F(n,x)=0
\label{eq:Legendre1}
\end{equation}
and their holonomic recurrence equation
\begin{equation}
(n+2)F(n+2,x)-(3 + 2 n) x  F(n+1,x)+(n+1)F(n,x)=0
\;,
\label{eq:Legendre2}
\end{equation}
together with the initial values
\begin{equation}
F(0,0)=1
\;,\quad
F(1,0)=0
\;,\quad
F'(0,0)=0
\;,\quad
F'(1,0)=1
\;.
\label{eq:Legendreinitial}
\end{equation}
Equations (\ref{eq:Legendre1})--(\ref{eq:Legendreinitial}) therefore build a
sufficient algebraic, even polynomial structure to represent the
functions $P_n(x)$ as we shall see now.

If we interpret the (partial) differentiations and shifts that occur 
as operators, and the representing system of holonomic equations as operator
equations, then these form a {\sl polynomial equations system} in a
noncommutative polynomial ring. For a continuous variable $x$ with
differential operator $D$ given by $DF(n,x)=F'(n,x)$,
the product rule implies \mbox{$D (x f)-x D f=f$,} and hence the
commutator rule $D x-xD=1$ is valid. Similarly for a discrete variable
$n$ with the (forward) shift operator $N$ given by $NF(n,x)=F(n+1,x)$,
we have $N (n F(n,x))-n N F(n,x)=$ $(n+1) F(n+1,x)-n F(n+1,x)=F(n+1,x)=
NF(n,x)$, and therefore the commutator rule $N n-n N=N$.
Similar rules are valid for all variables involved, 
whereas all other commutators vanish.

The transformation of a holonomic system given by mixed holonomic
differ\-ence-differential equations represents an elimination problem
in the noncommutative polynomial ring considered,
that can be solved by noncommutative Gr\"obner basis methods
(\cite{BW}, \cite{Gal}, \cite{Weispfenning}, \cite{Zei1}, \cite{Zei4},
\cite{Ta1}--\cite{Ta3}), \cite{Koe94z}).

Hence, we need the concept of a {\sl Gr\"obner basis}. If one applies 
Gau{\ss}'s
algorithm to a linear system, the variables are eliminated iteratively,
resulting in an equivalent system which is simpler in the sense that
it contains some equations which are free of some variables involved.
Note that connected with an application of Gau{\ss}'s algorithm is a
certain order of the variables.

The {\sl Buchberger algorithm}
is an elimination process, given a certain term order for the variables
(a variable order is no longer sufficient), 
with which a polynomial system (rather than a linear one)
is transformed, resulting in an equivalent system (i.e., constituting the
same ideal) for which the terms
that are largest with respect to the term order, are eliminated as far
as possible. Note that---in contrast to the linear case---the resulting 
equivalent system may contain more polynomials than the original one.
Such a rewritten system is called a Gr\"obner basis of the ideal
generated by the polynomial system given.
It turns out that Buchberger's algorithm can be extended to the
noncommutative case that we consider here \cite{Weispfenning} as long
as the rewrite process using the commutator does not increase the
variable order.

As an example, we consider $F(n,k)={{n}\choose{k}}$ in which case we have
the Pascal triangle relation $F(n+1,k+1)=F(n,k)+F(n,k+1)$, together with
the pure recurrence equation
$(n+1-k) F(n+1,k)-(n+1) F(n,k)=0$ with respect to $n$, say.
These equations read as $(KN-1-K) F(n,k)=0$,
and $((n+1-k)N-(n+1))F(n,k)=0$ in operator notation, $K$ denoting
the shift operator with respect to $k$. Therefore we have the polynomial system
\begin{equation}
KN-1-K
\quad\quad\mbox{and}\quad\quad
(n+1-k)N-(n+1)
\;.
\label{eq:leftideal1}
\end{equation}
The Gr\"obner basis of the left ideal generated by (\ref{eq:leftideal1})
with respect to the lexicographical term order $(k,n,K,N)$ is given by
\[
\Big\{(k+1) K +k-n,(n+1-k)N-(n+1),KN-1-K\Big\}
\;,
\]
i.e., the elimination process has generated the pure recurrence equation
\[
(k+1) F(n,k+1) +(k-n) F(n,k)=0
\]
with respect to $k$.

We used the \reduce\ implementation \cite{AM} for the noncommutative
\linebreak
Gr\"obner calculations of this article, but
I would like to mention that there is also 
\linebreak
a Maple package {\tt Mgfun}
written by Chyzak \cite{Chyzak} (to be obtained from
\linebreak
{\tt http://pauillac.inria.fr/algo/libraries/libraries.html\#Mgfun})
which can be used for this purpose.

As another example, we consider the Legendre polynomials. 
In operator notation the holonomic equations 
(\ref{eq:Legendre1})--(\ref{eq:Legendre2}) constitute the polynomials
\begin{equation}
(x^2-1)D^2+2xD -n(1+n)
\quad\mbox{and}\quad
(n+2)N^2-(3+2n)x N+(n+1)
\;.
\label{eq:Legendreoperator notation}
\end{equation}
The Gr\"obner basis of the left ideal generated by
(\ref{eq:Legendreoperator notation}) with respect to the lexicographical
term order $(D,N,n,x)$ is given by
\[
\Big\{
(x^2-1)D^2+2xD -n(1+n),
\]
\begin{equation}
(1+n)ND-(1+n)xD-(1+n)^2,
\label{eq:Legendre2a}
\end{equation}
\begin{equation}
(x^2-1)ND-(1+n)xN+(1+n),
\label{eq:Legendre3}
\end{equation}
\begin{equation}
(1+n)(x^2-1)D-(1+n)^2 N+x(1+n)^2,
\label{eq:Legendre4}
\end{equation}
\[
(n+2)N^2-(3+2n)x N+(n+1)
\Big\}
\;.
\]
After the calculation of the Gr\"obner basis, for better readability
I positioned the operators $D$ and $N$ back to the right, so that
the equations can be easily understood as operator equations, again.
By the term order chosen, the Gr\"obner basis contains those
equations for which the $D$-powers are eliminated as far as possible,
and (\ref{eq:Legendre2a})--(\ref{eq:Legendre4}) correspond to the
relations
\[
P_{n+1}'(x)=x\,P_{n}'(x)+(1+n)\,P_{n}(x)
\;,
\]
\[
(x^2-1)P_{n+1}'(x)=(1+n)\,(x P_{n+1}(x)-P_n(x))
\;,
\]
\begin{equation}
(x^2-1)P_{n}'(x)=(1+n)\,(P_{n+1}(x)-x P_n(x))
\; 
\label{eq:derivative rule}
\end{equation}
between the Legendre polynomials and their derivatives.

If we are interested in a relation between the Legendre polynomials
and their derivatives that is $x$-free (which is of importance
for example for spectral approximation, see \cite{CHQZ}), we choose
the term order $(x,D,N,n)$ to eliminate $x$ in the first place,
and obtain a different Gr\"obner basis containing the $x$-free polynomial
\[
-(n + 2)(n + 1)D-(2 n + 3) (n + 2) (n + 1) N+ (n + 2) (n + 1) N^2 D
\]
equivalent to the identity
\[
(2n+1)P_n(x)=P_{n+1}'(x)-P_{n-1}'(x)
\]
for the Legendre polynomials (see e.g.\ \cite{CHQZ}, formula (2.3.16)).

Here, we present the \reduce\ output for the above examples:

{\small
\begin{verbatim}
1: load ncpoly;

2: nc_setup({D,NN,n,x},{NN*n-n*NN=NN,D*x-x*D=1},left);

3: p1:=(x^2-1)*D^2+2*x*D-n*(1+n)$ % differential equation

4: p2:=(n+2)*NN^2-(3+2*n)*x*NN+(n+1)$ % recurrence equation

5: nc_groebner({p1,p2});

  2  2    2            2
{d *x  - d  - 2*d*x - n  - n,

                         2
 d*nn*n - d*n*x - d*x - n  - n,

       2
 d*nn*x  - d*nn - nn*n*x - 2*nn*x + n + 1,

      2            2           2    2
 d*n*x  - d*n + d*x  - d - nn*n  + n *x - x,

   2
 nn *n - 2*nn*n*x - nn*x + n + 1}

6: nc_setup({x,D,NN,n},{NN*n-n*NN=NN,D*x-x*D=1},left);

7: result:=nc_groebner({p1,p2});

            2  2            2    2
result := {x *d  + 2*x*d - d  - n  - n,

                        2                     2
           x*d*nn - d*nn *n + d*n + d + 2*nn*n  + 2*nn*n + nn,

                                   2
           x*d*n + x*d - d*nn*n + n  + 2*n + 1,

                               2
           2*x*nn*n + x*nn - nn *n - n - 1,

               2  2       2        2                       3         2
           d*nn *n  - d*nn *n - d*n  - 3*d*n - 2*d - 2*nn*n  - 3*nn*n  - nn*n}

8: nc_setup({n,x,NN,D},{NN*n-n*NN=NN,D*x-x*D=1},left);

9: nc_compact(part(result,5));

                                                    2
 - (2*n + 3)*(n + 2)*(n + 1)*nn + (n + 2)*(n + 1)*nn *d - (n + 2)*(n + 1)*d
\end{verbatim}
}\noindent
We see, therefore, that by the given procedure new relations
(between the binomial coefficients, and between the derivatives
of the Legendre polynomials) can be {\sl discovered}.
The generation of derivative rules like (\ref{eq:derivative rule}),
and the algorithmic work with them is described in \cite{Koe94z}.

\section{Holonomic Sums and Integrals}
\label{sec:Holonomic Sums and Integrals}

Analogously, with the method in the last section, holonomic
recurrence equations for {\sl holonomic sums} can be generated.
Note that the idea to use recurrence equations for the summand to
deduce a recurrence equation for the sum is originally due to
Sister Celine Fasenmyer (\cite{Fasenmyerdiss}--\cite{Fasenmyer}, 
see \cite{Rainville}, Chapter 14). Zeilberger \cite{Zei1}
brought this into a more general setting.

Consider for example
\[
s(n)=\sum_{k=0}^n F(n,k)=\sum_{k=0}^n {{n}\choose{k}} P_n(x)
\;,
\]
then by the product algorithm, we find the holonomic recurrence equations
\[
(n-k+1) F(n+1,k)-(1+n) F(n,k) = 0
\]
and
\[
(2+k)^2 F(n,k+2)-(3+2k)(n-k-1)x F(n,k+1)+(n-k)(n-k-1) F(n,k)=0
\]
for the summand $F(n,k)$. The Gr\"obner basis of the left ideal generated by
the corresponding polynomials
\[
(n-k+1) N-(1+n)
\quad\mbox{and}\quad
(2+k)^2 K^2-(3+2k)(n-k-1)x  K +(n-k)(n-k-1) 
\]
with respect to the lexicographical term order $(k,N,n,K)$ contains the
$k$-free polynomial
\begin{equation}
(2\!+\!n)^2 K^2 N^2 \!-\! K (2\! +\! n) (3 \!+\! 2 n) (K \!+\! x) N \!+ 
\!(1 \!+\! n) (2\! +\! n) (1\! +\! K^2 \! + \!2 K x)
\;,
\label{eq:kfree}
\end{equation}
which corresponds to a $k$-free recurrence equation for $F(n,k)$. We
use the order $(k,N,n,K)$ because then $k$-powers are eliminated as far as
possible (since we like to find a $k$-free recurrence), 
and $N$-powers come next in the elimination process (since the recurrence
equation obtained should be of lowest possible order).

Because all shifted sums
\[
s(n)=\sum_{k\in\Z} F(n,k)=\sum_{k\in\Z} F(n,k+1)=\sum_{k\in\Z} F(n,k+2)
\]
generate the same function $s(n)$, and since summing the $k$-free recurrence
equation is equivalent to setting $K=1$ in the corresponding operator equation
(check!), the substitution $K=1$ in (\ref{eq:kfree})
generates the valid holonomic recurrence equation
\[
(2+n) s(n+2) - (3 + 2 n) (1 + x) s(n+1) + 2 (1 + n) (1 + x) s(n)=0
\]
for $s(n)$.

In the general case, 
we search for a $k$-free recurrence equation contained in a Gr\"obner basis
of the corresponding left ideal with respect to a suitably chosen
weighted \cite{MNM} (or lexicographical $(k,N,n,K)$)
term order. For example, the elimination problems described
in \cite{Zei4} are automated by this procedure.

On the other hand, it turns out that in many cases the holonomic recurrence
equation derived is not of the 
lowest order. In the next section, we will discuss
how this problem can be resolved.

Note that by a similar technique, holonomic integrals can be treated
\cite{AZ}. To find a holonomic equation for
\[
I(y):=\int\limits_{a}^b
F(y,x)\,dx
\]
for holonomic $F(y,x)$ with respect to the discrete or continuous variable
$y$, calculate the Gr\"obner basis of the left ideal constituted by
the holonomic equations of $F(y,x)$ with respect to a suitably chosen
weighted or the lexicographical term order $(x,D_y,y,D_x)$.
We search for an $x$-free holonomic equation $\cal E$ contained in such 
a Gr\"obner basis. In case, that
$F(y,a)=F(y,b)\equiv 0$, and enough derivatives of $F(y,x)$ with respect
to $x$ vanish at $x=a$ and $x=b$, by partial integration it follows that
the holonomic equation valid for $I(y)$ is given by the substitution $D_x=0$ 
into $\cal E$ (see \cite{AZ}).

As an example, we consider
\[
I(n):=\int\limits_{-\infty}^\infty
e^{-x^2}\,H_n(x)\,dx
\;,
\]
$H_n(x)$ denoting the Hermite polynomials.
The method of \S~\ref{sec:Identification of Transcendental Functions}
yields the holonomic polynomials
\[
2\,(1 + n) + N^2 - 2\,x\,N
\quad\quad
\mbox{and}
\quad\quad
D^2 + 2\,(1 + n) + 2\,x\,D
\]
for the integrand.
Note that since $H_n(x)$ is an odd function for odd $n$, 
it is immediately clear that $I(n)=0$ in this case. 
However, what about even values of $n$?

The Gr\"obner basis of the corresponding left ideal contains the two
$x$-free polynomials
\[
N^2+N\,D
\quad\quad
\mbox{and}
\quad\quad
N\,n+n\,D+D
\]
so that setting $D=0$ we get for $I(n)$ the recurrence equation $I(n+1)=0$.
Indeed, this proves that $I(n)=0$ for $n\geq 1$.

As another example, we consider the Abramowitz functions (\cite{AS}, 27.5))
\[
A(n,y):=\int\limits_0^\infty x^n\,e^{-x^2-y/x}\,dx
\;.
\]
By the method in \S~\ref{sec:Identification of Transcendental Functions}
for the integrand $F(n,y,x)=x^n\,e^{-x^2-y/x}$
we get the three holonomic polynomials
\[
x-N
\;,
\quad\quad
-n\,x + x^2\,D_x + 2\,x^3 - y
\quad\quad
\mbox{and}
\quad\quad
1 + x\,D_y
\;.
\]
Using the term order $(x,D_y,y,D_x)$, the differential equation
\[
y\,A'''(n,y)-(n-1)\,A''(n,y)+2A(n,y)=0
\;,
\]
and using $(x,N,n,D)$, the recurrence equation 
\[
2A(n+3,y) - (n + 2)\,A(n+1,y) - y\,A(n,y)=0
\]
is automatically
generated by the given approach (compare \cite{AS}, (27.5.1), (27.5.3)).

Finally, we mention that similarly an identity like (\cite{AS}, (11.4.28))
\begin{equation}
\int\limits_0^\infty e^{-a^2\,x^2}\,x^{m-1}\,J_n(bx)\,dx
=
\frac{\Gamma\:(n/2+m/2)\,b^n}{2^{n+1}\,a^{n+m}\,\Gamma\:(n+1)}
\; _1 F_1\left.
\!\!
\left(
\begin{array}{c}
n/2+m/2
\\[1mm]
n+1
            \end{array}
\right| -\frac{b^2}{4a^2}\right)
\label{eq:Bessel}
\end{equation}
($ _1 F_1$ representing Kummer's confluent hypergeometric function)
for the Bessel function is proved by the calculation of the common
holonomic recurrence equation
\begin{eqnarray*}
0&=&
-\left (n+3\right )\left (n+m\right ){b}^{2}{I}(n)
\\&&
+2\,\left (n+
2\right )\left (4\,{a}^{2}{n}^{2}+16\,{a}^{2}n+12\,{a}^{2}-{b}^{2}m+{b
}^{2}\right ){I}(n+2)
\\&&
+\left (n+1\right )\left (n+4-m\right ){b}^{2}{I}(n+4)
\end{eqnarray*}
for the left and right hand sides of (\ref{eq:Bessel}). Note that
Zeilberger's algorithm is not directly applicable to the right hand side, but
the extended version of \cite{Koe94z} gives the result.

\section{Noncommutative Factorization and Holonomic Normal Form}
\label{sec:Noncommutative Factorization and Holonomic Normal Form}

Note that neither the sum and product algorithms of 
\S~\ref{sec:Identification of Transcendental Functions}, nor Zeilberger's
algorithm or its extension \cite{Koe94z},
nor the algorithms for holonomic sums and integrals
of \S~\ref{sec:Holonomic Sums and Integrals}
can guarantee to present the holonomic equation $\cal N$ of lowest order, and 
therefore the normal form searched for.

In \cite{MK}%
\footnote{Due to a severe bicycle accident of Herbert Melenk,
this paper is still unfinished.}
a Gr\"obner basis based factorization algorithm was introduced for polynomials
in noncommutative polynomial rings given by Lie bracket commutator rules.
This method is implemented in \cite{AM}. Given an expression $f$, and a
holonomic equation $\cal P$ of order $m$ of $f$,
one may find the normal form $\cal N$ of $f$
using this factorization algorithm
by generating the right factors of the noncommutative polynomial $p$
corresponding to $\cal P$, and checking if any of them, $\cal Q$, say,
(having order $l<m$, say)
and $m-l$ derivatives (shifts) of $\cal Q$ are satisfied by $f$
at a certain initial 
point. In the affirmative case, $\cal Q$ is compatible with $f$,
and corresponds to a valid holonomic equation for $f$.

To present some examples, we consider Zeilberger's algorithm first.
An example for which Zeilberger's algorithm does not generate the
holonomic recurrence equation of lowest order
is given by the sum (see e.g.\ \cite{PS})
\[
s_n:=\sum_{k=0}^n
(-1)^k {{n}\choose{k}}{{3k}\choose{n}}
\]
for which the holonomic equation
\begin{equation}
2\,\left (2\,n+3\right )\,s_{n+2}+3\,\left (5\,n+7\right )\,s_{n+1}
+9\,\left (n+1\right )\,s_n
=0
\label{eq:right factor}
\end{equation}
is generated. 
Note that there is an algorithm due to Petkov\u sek \cite{Pet} to find all
hypergeometric solutions of holonomic recurrence equations which
could be used as next step. However, we may also proceed as follows:
The corresponding noncommutative polynomial
$ 2(2n+3)N^2+3(5n+7)N+9(n+1)$
is factorized by implementation
\cite{AM} as
\[
2(2n+3)\,N^2+3(5n+7)\,N+9(n+1)=
((4n+6)\,N + 3(n+1))\,(N + 3)
\;.
\]
The right factor $N+3$ corresponds to the holonomic recurrence 
equation 
\begin{equation}
S_{n+1}+3S_n=0
\label{eq:right factor2}
\;,
\end{equation}
which, together with the initial value
$S_0=s_0=1$ uniquely defines a sequence $(S_n)_{n\in\N_0}$. Since
$S_1=-3$ turns out to be compatible with the given sum
\[
s_1=\sum_{k=0}^1
(-1)^k {{1}\choose{k}}{{3k}\choose{1}}
=-3
\;,
\]
and since (\ref{eq:right factor2}) implies (\ref{eq:right factor})
(right factor!), the sequence
$s_n$, which is the unique solution of (\ref{eq:right factor}) with $s_0=1$ and
$s_1=-3$, must equal $S_n$. From (\ref{eq:right factor2}), however,
the closed form $s_n=(-3)^n$ follows.

Similarly, for any particular $d\in\N$, $d\geq 3$, the identity
\[
\sum_{k=0}^n
(-1)^k {{n}\choose{k}}{{dk}\choose{n}}
=(-d)^n
\]
can be established, for whose left hand side
Zeilberger's algorithm generates a recurrence equation of order $d-1$
(see \cite{PS}).

Whereas Petkov\u sek's algorithm finds hypergeometric solutions
of holonomic recurrence equations as in the example, and
therefore not only verifies identities, but {\sl generates}
closed-form results, our approach is more general in the following sense.
Factorizations with polynomial coefficients
of ordinary holonomic differential equations 
(see \cite{Beke1}, \cite{Schwarz} for other methods)
as well as of any mixed holonomic difference-differential equation
can be calculated.

We give an example of that type for the sum algorithm:
Consider the difference of successive
Gegenbauer polynomials $h(x)=C_{n+1}^{(-1/2)}(x)-C_n^{(-1/2)}(x)$ 
that were used in \cite{KS:Koebe}.
Here the summand $f(x):=C_n^{(-1/2)}(x)$ satisfies the holonomic
equation
\[
( {x^2}-1 ) \,f''(x)+( n - {n^2} ) \,f(x) = 0
\;,
\]
and the sum algorithm yields the fourth order equation
\[
    {{( x^2-1 ) }^2}\,h''''(x) + 4\,x\,( x^2-1 ) \,h'''(x)
-
    2\,(n^2-1) ( x^2-1 ) \,h''(x) + {n^2}\,( n^2-1 ) \,h(x) 
=0
\]
for $h(x)$. The implementation \cite{AM} finds (besides others)
the noncommutative factorization 
\[
\Big((x^2-1)\,D^2+(1+x)\,D-n^2\Big)\,\Big((x^2-1)\,D^2-(1+x)\,D+(1-n^2)\Big)
\]
of the corresponding noncommutative polynomial, whose right factor
\[
(x^2-1)\,D^2-(1+x)\,D+(1-n^2)
\]
turns out to be compatible with the given function $h(x)$. That is, the
corresponding differential equation and two derivatives thereof
are satisfied by $h(x)$, at $x=1$.
Therefore the holonomic normal
form of $h(x)$ is the corresponding differential equation 
\[
  ( 1 - {n^2} ) \,h(x) - \left( 1 + x \right) \,h'(x) + 
    ( {x^2}-1 ) \,h''(x) = 0
\]
that was a tool in \cite{KS:Koebe}. This result can also be obtained by the
method given in \cite{Koe92}--\cite{Koe93}.

To evaluate the integrals
\[
I_n:=\int\limits_{-\infty}^\infty
x^n\,e^{-x^2}\,H_n(x)\,dx
\;,
\]
we may deduce the holonomic system
\[
N^2-2x^2N+2(1+n)x^2
\]
and
\[
x^2\,D^2 + 2x(x^2-n)\,D+(n+n^2+2x^2)
\]
of the integrand.  The Gr\"obner basis of this system
with respect to the weighted lexicographical order with weights
$(3,1,0,0)$ for $(x,N,n,D)$ (i.e.\ the term
$x$ is considered larger than $N^3$,
whereas $x$ is smaller than $N^4$, and any power of $n$ and $D$ is smaller than
$x$ and $N$) contains an $x$-free polynomial, which when
evaluated at $D=0$ yields
\begin{eqnarray}
P(n,N)&=&(n + 5) (n + 4) (n + 3) N^3
-(3 n + 7) (n + 5) (n + 4) (n + 3) N^2
\nonumber
\\&&+(3 n + 5) (n + 5) (n + 4) (n + 3) (n + 2) N
\label{eq:factorization}
\\&&-
 (n + 5) (n + 4) (n + 3) (n + 2) (n + 1)^2
\nonumber
\end{eqnarray}
corresponding to a recurrence equation of order three. 

On the other hand, $P(n,N)$ obviously has the trivial 
(commutative) factorization
\[
P(n,N)=
(n + 5) (n + 4) (n + 3)
\Big(N^3-(3 n + 7)N^2+(3 n + 5) (n + 2) N- (n + 2) (n + 1)^2\Big)
\]
and the remaining right factor can be represented as
\[
N^3-(3 n + 7)N^2+(3 n + 5) (n + 2) N- (n + 2) (n + 1)^2
=
(N-n-2)(N-n-1)(N-n-1)
\]
(note that \cite{AM} finds four different right factors).
This leads to the valid recurrence equation $I_{n+1}=(n+1)I_n$
that together with the initial value 
\[
I_0=
\int\limits_{-\infty}^\infty
e^{-x^2}\,dx
=\sqrt \pi
\]
gives finally $I_n=\sqrt \pi\,n!$.

\section*{Acknowledgement}
I would like to thank Prof.\ Peter Deuflhard for his support and
encouragement. I'd also like to thank
Herbert Melenk for his advice on Gr\"obner bases,
and for his excellent \reduce\ implementation \cite{AM}.
Hopefully, he will have recovered soon!

\end{document}